\documentclass[12pt,leqno]{article}
\usepackage{amssymb, amscd, amsmath, amsthm}

\def\beq{\begin{equation}}
\def\eeq{\end{equation}}

\theoremstyle{definition}

\allowdisplaybreaks

\theoremstyle{theorem}

\newtheorem{theorem}{Theorem}

\newtheorem{corollary}{Corollary}

\theoremstyle{definition}
\theoremstyle{remark}
\newtheorem{remark}{Remark}

\begin{document}

\begin{center} 

{\large{On the frequency of small gaps between primes}}

\medskip 
\centerline{by}

\medskip 
{\large{\bf \'Akos Magyar}} 

\medskip 
and 

\medskip 
{\large {\bf J\'anos Pintz}}

\medskip 
(HUN-REN Alfr\'ed R\'enyi Institute of Mathematics 

and the University of Georgia)

\end{center}

\bigskip 
\begin{abstract}
In a recent work Friedlander studied the problem of how large
consecutive prime gaps should be in order that the sum of the
reciprocals should be divergent.
Supposing a very deep Hypothesis, a generalization of the Hardy--Littlewood prime $k$-tuple conjecture, he gave an almost precise answer for it.
In the present work we give an unconditional answer for a much weaker form of the same problem.

\footnote{The first author was supported by Grant {\#}854813 by the Simons Foundation.
The second author was supported by the Hungarian National Research Developments and Innovation Office, NKFIH, 147153 and KKP133819, Excellence 151341}

\footnote{Key words and phrases: Brun's theorem on twin primes, reciprocal of consecutive
primes, small gaps between primes, frequency of small prime gaps.}

\footnote{2020 {Mathematical Subject Classification: Primary 11N05, Secondary 11N36.}}

\end{abstract}

\setcounter{equation}{0}
\numberwithin{equation}{section}

\setcounter{section}{1}
\noindent
John Friedlander investigated in his recent work entitled ``When Euler met Brun''
\cite{Fri2024} that how small one can choose $\lambda(p)$ in order that
\beq
\label{eq:1.1}
\sum_{p' - p \leq \lambda(p)\log p} \frac1{p} = \infty
\eeq
should still hold, where $p'$ denotes the prime following $p$ ($p$ runs over all primes).

\vfil\goodbreak
The title is motivated by the famous result of Euler [{Eul1737}, Theorem~7]
\beq
\label{eq:1.2}
\sum_p \frac1{p} = \infty
\eeq
and that of Brun \cite{Bru1915} that restricting the summation for the twin primes the corresponding sum is already convergent.
Brun's method immediately yields that for every fixed $K$ one has also
\beq
\label{eq:1.3}
\sum_{p' - p \leq K} \frac1{p} < \infty.
\eeq

In order to determine the threshold of the function $\lambda(p)$ Friedlander introduces the function
\beq
\label{eq:1.4}
\log_k x = \log\log_{k - 1} x, \ \ \ \log_1 x = \log x
\eeq
and what he calls the logorial function
\beq
\label{eq:1.5}
\text{\rm Log}_k x = \prod_{2 \leq j \leq k} \log_j x.
\eeq

Using the above notation he proved [{Fri2024}, (3.8)] without any condition that with
$y(p) = \lambda(p) \log p$
\beq
\label{eq:1.6}
\sum_{\substack{p\\ p' - p \leq y(p)}} \frac1{p} < \infty \ \text{ if } \ \lambda(p) = 1\!\bigm/\!(\text{\rm Log}_k p)(\log_k p)^\varepsilon
\eeq
for any $\varepsilon > 0$.
If he supposes the generalization of a uniform version of Hardy--Littlewood's prime $k$-tuple conjecture (cf. Hypothesis G below), then using a theorem of Gallagher \cite{Gal1976} (and a corollary of that by K. Soundararajanan [{Sou2007}, Exercise 1.3], see also Goldston and Ledoan's work \cite{GL2012/13}) Friedlander derives the conditional result
\beq
\label{eq:1.7}
\sum_{\substack{p \leq x\\
p' - p \leq y(p)}} \frac1{p} \sim \log_{k + 1} x \to \infty \ \text{ for } \ \lambda(p) = 1\!\bigm/\!\text{\rm Log}_k p
\eeq
for any $k \geq 2$.

In order to formulate the Hypothesis G, let $\mathcal H = \{h_1, \ldots, h_k\}$ denote a set of distinct non-negative integers and let $\nu_{\mathcal H}(p)$ denote the number of residue classes modulo $p$ occupied by the $h_i$.
We call $\mathcal H$ \emph{admissible} if $\nu_{\mathcal H}(p) < p$ for all primes $p$.
The singular series $\mathfrak S(\mathcal H)$ is defined by
\beq
\label{eq:1.8}
\mathfrak S(\mathcal H) = \prod_p \left(1 - \frac{\nu_{\mathcal H}(p)}{p}\right) \left(1 - \frac1{p}\right)^{-k}.
\eeq
This is positive iff the $k$-tuple $\mathcal H$ is admissible.

\smallskip
\noindent
{\bf Hypothesis G.}
{\it Let $\pi(x, \mathcal H)$ denote the number of integers $n \leq x$ such that $n + h_j$ is prime for all $h_j \in \mathcal H$.
Then
\beq
\label{eq:1.9}
\pi(x, \mathcal H) \sim \mathfrak S(\mathcal H) \frac{x}{(\log x)^k} \ \ \text{ as } \ x \to \infty
\eeq
uniformly for all $h_1,\ldots, h_k \leq \lambda \log x$, where $\lambda > 0$ is fixed.}

\smallskip
\noindent
{\bf Theorem} (Gallagher){\bf.}
{\it Let $P_k(h, x)$ denote the number of integers $n \leq x$ such that $(n, n + h]$ contains exactly $k$ primes.
Then, under Hypothesis~G we have for $h = \lambda \log x$, $\lambda > 0$ fixed
\beq
\label{eq:1.10}
P_k(h, x) \sim x e^{-\lambda} \frac{\lambda^k}{k!} \ \ \text{ as } \ x \to \infty.
\eeq
}

\smallskip
Soundararajanan (and also Goldston and Ledoan) derived from this the following

\smallskip
\noindent
{\bf Corollary.}
{\it Under Hypothesis G, for fixed $\lambda > 0$ we have
\beq
\label{eq:1.11}
\frac1{\pi(x)} \# \left\{ p \leq x; \frac{p' - p}{\log p} \leq \lambda\right\} \sim \int\limits_0^\lambda e^{-t} dt = 1 - e^{-\lambda}.
\eeq
}

\smallskip
Friedlander supposed in his work a generalization of the above Corollary when $\lambda(x)$ is allowed to tend to zero, if this convergence is sufficiently slow.

\smallskip
\noindent
{\bf Hypothesis S.} {\it If $\lambda(u) \to 0$ as $u \to \infty$, but
\beq
\label{eq:1.12}
\lambda(u) \gg \frac1{(\log\log u)^2},
\eeq
then \eqref{eq:1.11} remains valid with $\lambda(x)$ instea of $\lambda$.}

\smallskip
He also notes that for $\lambda \to 0$ we have $1 - e^{-\lambda} \sim \lambda$.
The unconditional result \eqref{eq:1.6} and the conditional \eqref{eq:1.7}
(see (3.8)--(3.9) of \cite{Fri2024}) nearly exactly determine under Hypothesis~H the threshold function $\lambda(p)$, therefore $y(p)$, for which the sum
\beq
\label{eq:1.13}
\sum_{\substack{p\\
p' - p \leq y(p)}} \frac1{p}
\eeq
is convergent, or divergent, respectively.

As quoted in \cite{Fri2024}, an unconditional lower bound for the number of prime gaps smaller than $\lambda(p)\log p$ was proved in \cite{GPY2013} when $\lambda(p) = \lambda > 0$ is an arbitrary fixed constant.
Theorem~1 of \cite{GPY2013} showed
\beq
\label{eq:1.14}
\# \{p \leq x; p' - p \leq \lambda \log p \} \gg_\lambda \pi(x)
\eeq
which shows that a positive proportion of prime gaps is smaller than $\lambda \log p$ for any fixed $\lambda > 0$.
On the other hand it was also proved in [{GPY2013}, Theorem 3] that Selberg's upper bound sieve coupled with Gallagher's theorem \cite{Gal1976}
\beq
\label{eq:1.15}
\sum_{\substack{\mathcal H \subset [1, h]\\
|\mathcal H| = k}} \mathfrak S(\mathcal H) \sim h^k \ \ \text{ for fixed } \ k \ \text{ as } \ h \to \infty,
\eeq
together yield for any $h > 2$
\beq
\label{eq:1.16}
\# \{p \leq x, p' - p \leq h\} \ll  \min \left(\frac{h}{\log x}, 1\right) \pi(x).
\eeq
In particular, if $h = o(\log x)$
\beq
\label{eq:1.17}
\# \{p \leq x; p' - p \leq h\} = o(\pi(x)).
\eeq

The goal of our present work is to extend \eqref{eq:1.14} for the case of $\frac{c}{\log p}
\leq \lambda(p) \ll  1$ with an unconditional explicit lower bound on the RHS of \eqref{eq:1.14} in dependence of the function $\lambda(p)$.
Theorem 2 of \cite{GPY2013} contained already such an explicit estimate for $\lambda > 0$ fixed but the lower bound
\beq
\label{eq:1.18}
\exp(-c\lambda^{-6})\pi(x)
\eeq
was rather weak.

In the present work we will show the following

\begin{theorem}
\label{th:1}
If $N$ is sufficiently large, that is $N > C_1$ (explicitly calculable),
$\eta = \eta(N) > C_0/\log N$ ($C_0$ explicitly computable positive constant), then
\beq
\label{eq:1.19}
\sum_{\substack{p \sim N\\
p' - p \leq \eta \log N}} 1 > c_0 \eta^{1879} \pi(N)
\eeq
with an explicitly calculable $c_0 > 0$, where $n \sim N$ means $N \leq n < 2N$.
\end{theorem}

This result gives an improved generalization of the result \eqref{eq:1.14} of \cite{GPY2013} about the positive proportion of small gaps between primes (where small means of size $\leq \varepsilon \log p$ for a given fixed $\varepsilon > 0$) and Zhang's result \cite{Zha2014} about the existence of infinitely many bounded gaps between primes.
\eqref{eq:1.19} is namely valid for all gap sizes between $C_0$ and $c'(\log N)$.

Using the notation \eqref{eq:1.5} we easily obtain a weaker unconditional version of \eqref{eq:1.7}.

\begin{corollary}
\label{cor:1}
Let $y(p) = \lambda(p) \log p$.
Then we have for any $m$
\beq
\label{eq:1.20}
\sum_{\substack{p\\
p' - p \leq \lambda(p)\log p}} \frac1{p} = \infty \ \ \text{ if }\ \lambda(p) \geq \frac{c'}{(\text{\rm Log}_m p)^{1/1879}}.
\eeq
\end{corollary}

\begin{corollary}
\label{cor:2}
For a $H > C_2$ (explicit absolute constant) we have below $N$
at least $C_3 H^{1879}N/(\log N)^{1880}$ prime gaps of size at most $H$.
\end{corollary}

\begin{remark}
\label{rem:1}
The exponent 1879 in Theorem \ref{th:1} (and therefore in Corollaries \ref{cor:1}
and \ref{cor:2}) can probably be reduced to a size around 50 if we use instead of
the method of Zhang and Polymath 8a the method of Maynard \cite{May2015} and Polymath 8B
\cite{Pol2014b}, \cite{Pol2015}.
We plan to return to this problem later.
\end{remark}

\begin{remark}
\label{rem:2}
The present proof does not produce the result above $N > C_1$
explicitly calculable positive constant.
However one can easily prove that a small modification, deleting the multiples of the maximal prime
divisor of a possible  (unique) exceptional modulus from the sieving set
of $d$'s changes the size of the critical quantities $S_0(\mathcal H)$ and
$S_1(\mathcal H, h^*)$ in \eqref{eq:3.3}--\eqref{eq:3.4} by a factor $(1+o(1))$, so the resulting
Bombieri--Vinogradov type Theorems turn to be effective and the proof
will stay valid.
\end{remark}

\section{Main ideas of the proof. Notation}
\label{sec:2}

In this section we outline the main ideas of the proof and introduce the notation used.
Our present work will be based on ideas and results of the two earlier works \cite{Pin2010} and \cite{GPY2013}, but the notation will just in part follow theirs.
Concerning the improvement the most crucial tool will be Zhang's result \cite{Zha2014} in the improved form by the Polymath 8 project \cite{Pol2014a} about equidistribution of primes in arithmetic progressions which led to the proof of the existence of infinitely many bounded gaps between primes by Zhang \cite{Zha2014}.

The starting point was to use a form of Selberg's upper bound sieve to detect primes in ``admissible'' $k$-tuples $\mathcal H = \{h_1, h_2, \ldots, h_k\}$, more exactly in the sequence $n + \mathcal H$ where $n$ runs through integers in an interval of type $[N, 2N)$ with a sufficiently large integer $N$, which will be abbreviated by $n \sim N$.
In several applications of the method the parameter $k$ was chosen bounded as $N \to \infty$.
A hopeless final goal could be (at least some good approximation) of Hardy--Littlewood's prime $k$-tuple conjecture, i.e., \eqref{eq:1.9} for bounded $h_k$.

The tool, investigated in several works of Goldston and
Y{\i}ld{\i}r{\i}m (\cite{GY2003}, \cite{GY2007a}, \cite{GY2007b}) was to consider the primes with the weights of Selberg's upper bound sieve,
\beq
\label{eq:2.1}
\Lambda_R(n; \mathcal H, k) := \sum_{\substack{d\\ d \mid P_{\mathcal H}(n)}} \mu(d) \log^k \frac{R}{d},
\eeq
where
\beq
\label{eq:2.2}
P_{\mathcal H}(n) := \prod_{i = 1}^k (n + h_i), \ \ \ R\leq N^{1/2} \mathcal L^{-c}, \ \ \mathcal L = \log N.
\eeq

This first attempt (in an unpublished manuscript) already improved the 15-year-old record of Helmut Maier \cite{Mai1988} about
\beq
\label{eq:2.3}
\Delta = \liminf_{n \to \infty} \frac{p_{n + 1} - p_n}{\log p_n} \leq 0.2484\ldots,
\eeq
but failed to prove $\Delta = 0$, that is, the existence of an infinite gap sequence of order $o(\log p)$ (which we will call informally \emph{small gaps}).

It was a little later when Goldston, Pintz and Y{\i}ld{\i}r{\i}m \cite{GPY2009} used instead of \eqref{eq:2.1} a $k + \ell$-dimensional sieve with $\ell > 0$,
\beq
\label{eq:2.4}
\Lambda_R(n; \mathcal H, k + \ell) = \sum_{\substack{d\\ d \mid P_{\mathcal H}(n)}} \mu(d) \log^{k + \ell} \frac{R}{d}
\eeq
for the original $k$-dimensional problem which led to
\beq
\label{eq:2.5}
\Delta = 0.
\eeq

We remark here that in order to reach $\Delta = 0$ one needed to choose $k \to \infty$, $\ell \to \infty$, $\ell = o(k)$, but in order to show $\Delta < \varepsilon$, $k$ (and $\ell$) could stay below some $C(\varepsilon)$.
Using the very usual smooth weights
\beq
\label{eq:2.6}
\theta(m) = \begin{cases}
\log p &\text{ if }m = p \text{ is prime,}\\
0 &\text{ otherwise,}
\end{cases}
\eeq
we investigate in the present work the two crucial sums
(similarly to \cite{GPY2009}, \cite{GPY2010}, \cite{Pin2010}, \cite{GPY2013}, \cite{Zha2014}, \cite{Pol2014a})
\beq
\label{eq:2.7}
S_0(\mathcal H) := \sum_{n \sim N} \Lambda_R^2 (n; \mathcal H, k + \ell) \ \ \ (1 \leq \ell \leq k)
\eeq
and
\beq
\label{eq:2.8}
S_1(\mathcal H, h^*) := \sum_{n \sim N} \theta (n + h^*) \Lambda_R^2 (n; \mathcal H, k + \ell),
\ \ h^* \ll  \mathcal L := \log N.
\eeq

Here, similarly to \cite{Zha2014} and \cite{Pol2014a}, and somewhat similarly to \cite{Pin2010}, $k$ and $\ell$ are bounded integers as $N \to \infty$.
In the present work we will choose them as
\beq
\label{eq:2.9}
k = 1880 \ \ \ \text{ and } \ \ \ \ell = 21.
\eeq

On the other hand, similarly to \cite{GPY2009}, \cite{Pin2010} and \cite{GPY2013}, but in contrast to \cite{Zha2014} and \cite{Pol2014a} we will deal with $k$-tuples $\mathcal H$ satisfying
\beq
\label{eq:2.10}
\mathcal H \subseteq [1, h] \ \ \text{ with } \ \ h \ll  \mathcal L
\eeq
(while \cite{Zha2014} and \cite{Pol2014a} work with a bounded $h$).

We mentioned already the first main idea, to weight primes with the sieve weights \eqref{eq:2.1} which are optimal in some problems.
The modified weights \eqref{eq:2.4} turned out to be successful in showing the Small Gap Conjecture, i.e., $\Delta = 0$ in \eqref{eq:2.5} (cf.\ \cite{GPY2009}).
The third main idea, observed in \cite{Pin2010} was that the sieve weights in \eqref{eq:2.4} are concentrated on almost prime values of $P_{\mathcal H}(n)$ (see \eqref{eq:2.2}), which expressed in an exact form appears in Lemmas 1 and 2, that is, formulae \eqref{eq:3.1} and \eqref{eq:3.2} of the present work (taken over from Lemmas 3 and 4 from \cite{Pin2010}), where, with the notation
\beq
\label{eq:2.11}
P(m) = \prod_{p \leq m} p
\eeq
we show that the contribution of integers $n$ with
\beq
\label{eq:2.12}
\bigl(P_{\mathcal H}(n), P(R^\delta)\bigr) > 1
\eeq
is negligible to the sums on the LHS of \eqref{eq:3.1}--\eqref{eq:3.2} for any $\delta > \widetilde c > 0$.

The fourth main idea is that in the complementary case of almost primes $k$-tuples, i.e., if
\beq
\label{eq:2.13}
\bigl(P_{\mathcal H}(n), P(R^\delta)\bigr) = 1
\eeq
we have
\beq
\label{eq:2.14}
\Lambda_R(n; \mathcal H, k + \ell) \ll _{k, \ell, \delta} (\log R)^{k + \ell} \ \ \text{ if } \ R \gg N^{c_0}.
\eeq
If we choose, say, $c_0 = 1/5$ and $\delta$ a small positive absolute constant, say $10^{-10}$, then together with the choices of $k$ and $\ell$ in \eqref{eq:2.9} the implied constant (in \eqref{eq:2.14}) will be an absolute one.

Concerning bounded gaps between primes, a conditional result was proved already in the first work of Goldston, Pintz and Y{\i}ld{\i}r{\i}m \cite{GPY2009} leading to the unconditional result $\Delta = 0$.
It was shown there that if the Bombieri--Vinogradov theorem can be extended to $q \leq x^\vartheta$,
$\vartheta > 1/2$, we obtain infinitely many bounded prime gaps.

The observation of Motohashi and Pintz \cite{MP2008} (timely before 2010 and 2013, the time of the mentioned 3\textsuperscript{rd} and 4\textsuperscript{th} main ideas but to be used for results about prime gaps just later by \cite{Zha2014} and \cite{Pol2014a}) was that in order to obtain the Bounded Gap Conjecture it is sufficient to show the extension of the Bombieri--Vinogradov theorem for $x^\delta$ smooth values of the moduli (i.e., when $q$ has no prime factor larger than $x^\delta$).

Although this is still open until today, Zhang \cite{Zha2014} succeeded to show that this holds if the residues $a_q$ are restricted to the solution of an arbitrary but fixed polynomial congruence
\beq
\label{eq:2.15}
\sum_{\substack{q \leq x^{\vartheta - \varepsilon}\\
q\ x^\delta\text{\rm -smooth, } \mu(q) \neq 0}}
\left|\psi (x, q, a_q) - \frac{x}{\varphi(q)}\right| \ll  \frac{x}{(\log x)^A}.
\eeq
In his case the polynomial congruence was $\prod\limits_{1 \leq j \leq k, j \neq i} (n + h_j - h_i) \equiv 0 (\text{\rm mod}\,q)$.
He proved it with the value $\vartheta = 1/2 + 1/584$.
During the Polymath 8A project \cite{Pol2014a} this was shown in a more general and sharper form as

\smallskip
\noindent
{\bf Theorem 0} (Polymath 8A). {\it Let $\vartheta = 1/2 + 7/300$.
Let $\varepsilon > 0$ and $A \geq 1$ be fixed real numbers.
For all primes $p$, let $a_p$ be a fixed invertible residue class modulo $p$ and for $q \geq 1$ squarefree, denote by $a_q$ the unique invertible residue class modulo $q$ such that $a_q \equiv a_p (\text{\rm mod }p)$ for all primes $p$ dividing $q$.
Then \eqref{eq:2.15} holds with some $\delta > 0$, depending only on $\varepsilon$, and the constant implied by the $\ll $ symbol depends on $A$, $\varepsilon$ and $\delta$, but it is independent of the residue classes $(a_p)$.}

\smallskip
This is the most important result that we will use in the proof of Theorem \ref{th:1}.
We do not need the corresponding Bombieri--Vinogradov theorem for other residue classes
as shown in \cite{MP2008} and \cite{Zha2014}.

We allowed us the freedom to quote \eqref{eq:3.2} and \eqref{eq:3.4} (so as a consequence \eqref{eq:3.6}) with the condition $R \leq N^{(\vartheta - \varepsilon)/(2 + \delta)}$ as Lemmas 4 and 2 of \cite{Pin2010}, although the condition there was the validity of the full extension of the Bombieri--Vinogradov theorem.
However, as explained above, using instead \eqref{eq:2.15} and the above remarks, our Lemmas 2, 4 (and so 6) remain valid with $\vartheta = 1/2 + 7/300$ if the moduli are $\delta$-smooth and some polynomial congruence is satisfied for the relevant residues.

\section{Weighted sums of prime $k$-tuples}
\label{sec:3}

The structure of the proof will follow that of \cite{GPY2013} but the use of the strong estimate of Theorem~1 of \cite{Pol2014a} (see our Theorem 0 in Section \ref{sec:2}) will be decisive in reaching our goal.
The crucial Propositions 1 and 2 of \cite{GPY2013} are based on Lemmas 1, 2, 4 and 5 of \cite{Pin2010}.
However, we can use Lemma 5 now in a stronger form due to the mentioned result of \cite{Pol2014a}
(see the following Lemma 2), while we quote Lemma 4 of \cite{Pin2010} simply as

\smallskip
\noindent
{\bf Lemma 1.} {\it Let $N^{c_1} < R \leq N^{1/(2 + \delta)}(\log N)^{-C_2}$.
Then we have
\beq
\label{eq:3.1}
\sum_{\substack{n \sim N\\
(P_\mathcal H(n), P(R^\delta)) > 1}} \Lambda_R^2(n; \mathcal H, k + \ell) \ll  \delta \sum_{n \sim N} \Lambda_R^2 (n; \mathcal H, k + \ell).
\eeq
}

\smallskip
\noindent
{\bf Lemma 2.} {\it Let $N^{c_1} \leq R \leq N^{(\vartheta - \varepsilon)/(2 + \delta)}(\log N)^{-c_2}$, $0 < \delta < c_0$, $\varepsilon > 0$.
Let $h^* \leq h$ and $\kappa = 1$ if $h^* \in \mathcal H$, $\kappa = 0$ if $h^* \notin \mathcal H$.
Then
\begin{align}
\label{eq:3.2}
&\sum_{\substack{n \sim N\\
(P_{\mathcal H}(n), P(R^\delta)) > 1}} \theta(n + h^*) \Lambda_R^2(n,\mathcal H, k + \ell)\\
&\ll \delta \sum_{n \sim N} \theta(n + h^*) \Lambda_R^2(n, \mathcal H, k + \ell)\nonumber\\
&\quad  + O\bigl(N(\log_2 N)^c\bigl((\log N)^{k + 2\ell + \kappa - 1} + (\log N)^{k + \ell - 1/2}\bigr)\bigr).
\nonumber
\end{align}
}

\smallskip
The crucial quantities on the right-hand side of \eqref{eq:3.1} and \eqref{eq:3.2} were evaluated already in the works \cite{GPY2009} and \cite{GPY2010}.
We quote them now in a simplified form of Lemmas 1 and 2 of \cite{Pin2010} as

\smallskip
\noindent
{\bf Lemma 3.} {\it Let $N^{c_1} \leq R \leq N^{1/2} \log N^{-c_2}$.
Then
\beq
\label{eq:3.3}
S_0(\mathcal H) := \sum_{n \sim N} \Lambda_R^2 (n; \mathcal H, k + \ell) =
(1 + o(1)) \frac{\mathfrak S(\mathcal H)}{(k + 2\ell)!} {2\ell \choose \ell} N(\log R)^{k + 2\ell}.
\eeq
}

\smallskip
\noindent
{\bf Lemma 4.} {\it Let $N^{c_1} \leq R \leq N^{(\vartheta - \varepsilon)/2}$.
Let $h^* \leq h$ and $\kappa = 1$ if $h^* \in \mathcal H$, $\kappa = 0$ if $h^* \notin \mathcal H$.
Then
\begin{align}
\label{eq:3.4}
&S_1(\mathcal H, h^*)
:= \sum_{n \sim N} \theta(n + h^*) \Lambda_R^2 (n; \mathcal H, k + \ell) \\
&= \frac{\mathfrak S(\mathcal H \cup \{h^*\})}{(k + 2\ell + \kappa)!} {2(\ell + \kappa)\choose \ell + \kappa} N(\log R)^{k + 2\ell + \kappa}\\
&\qquad  +  O\bigl(N(\log N)^{k + 2\ell + \kappa - 1} (\log_2 N)^c\bigr).\nonumber
\end{align}
}

\smallskip
We repeat here that our $\vartheta$ is defined in a different way from the $\vartheta$ in \cite{Pin2010} (which is the distribution level of the primes), but the proofs of Lemmas 2 and 4 are the same as that of Lemmas 4 and 2 of \cite{Pin2010} since the equidistribution is required only in the residue classes $a(q)$ which are solutions of $P_{\mathcal H}(n) \equiv 0(q)$.

Lemmas 1--4 immediately imply Lemmas 5 and 6.

\smallskip
\noindent
{\bf Lemma 5.} {\it Let $N^{c_1} \leq R \leq N^{1/(2 + \delta)}\log N^{-c_2}$.
Then we have
\begin{align}
\label{eq:3.5}
S_0(\mathcal H, \delta) &:= \sum_{\substack{n \sim N\\
(P_{\mathcal H}(n), P(R^\delta)) = 1}} \Lambda_R^2 (n; \mathcal H, k + \ell)\\
&\ = (1 + O(\delta)) \frac{\mathfrak S(\mathcal H)}{(k + 2\ell)!} {2\ell\choose \ell} N(\log R)^{k + 2\ell}.\nonumber
\end{align}
}

\smallskip
In case of $\mathfrak S(\mathcal H) = 0$ the RHS is $o\bigl(N(\log R)^{k + 2\ell}\bigr)$.

\smallskip
\noindent
{\bf Lemma 6.} {\it Let $N^{c_1} \leq R \leq N^{(\vartheta - \varepsilon)/(2 + \delta)}(\log N)^{-c_2}$.
Let $h^* \leq h$ and $\kappa = 1$ if $h^* \in \mathcal H$, $\kappa = 0$ if $h^* \notin \mathcal H$.
Then
\begin{align}
\label{eq:3.6}
S_1(\mathcal H, \delta) &:= \sum_{\substack{n \sim N\\
(P_\mathcal H(n), P(R^\delta)) = 1}}
\theta(n + h^*) \Lambda_R^2 (n; \mathcal H, k + \ell)\\
&\ = (1 + O(\delta)) \frac{\mathfrak S(\mathcal H \cup  \{h^*\})}{(k + 2\ell + \kappa)!} {2(\ell + \kappa)\choose \ell + \kappa} N(\log R)^{k + 2\ell + \kappa}.\nonumber
\end{align}
In case of $\mathfrak S(\mathcal H \cup h^*) = 0$ the RHS is $o\bigl(N(\log R)^{k + 2\ell + \kappa} \bigr)$.
}

\section{Proof of Theorem 1}
\label{sec:4}

Let us denote
\beq
\label{eq:4.1}
Q(N, h) := \sum_{\substack{n \sim N\\
\pi(n + h) - \pi(n) > 1}}\!\!\! 1, \qquad \Theta(n, h) := \sum_{h^* = 1}^h \theta(n + h^*).
\eeq
If $n$ is counted by the sum, then we have a prime $p_j$ with $n < p_j < p_{j + 1} \leq n + h$.
Thus $p_{j + 1} - p_j < h$ and $p_{j + 1} - h \leq n < p_j$ so that to every such gap we have at most $\lfloor h \rfloor$ such integers $n$.
Therefore we obtain from the PNT with error term used at the ends of such gaps:
\beq
\label{eq:4.2}
Q(N, h) \leq h \sum_{\substack{N < p_j \leq 2N\\ p_{j + 1} - p_j \leq h}} 1 + O\bigl(N e^{-c\sqrt{\log N}}\bigr).
\eeq
We will investigate the quantity
\beq
\label{eq:4.3}
S_2(h, \delta) := \frac1{Nh^k(\log R)^{k + 2\ell}} \sum_{n \sim N} (\Theta(n, h) - \log 3N) \Bigl(
{\sum_{\mathcal H}}^* \Lambda_R^2 (n, \mathcal H, k + \ell)\Bigr)
\eeq
where
\beq
\label{eq:4.4}
{\sum_{\mathcal H}}^* := \sum_{\substack{\mathcal H\subseteq [1, h], |\mathcal H| = k\\
\mathcal H \text{ admissible}\\
(P_{\mathcal H}(n), P(R^\delta)) = 1}}.
\eeq

From Lemmas 5, 6 and \eqref{eq:1.15} we obtain a lower bound for $S_2(h, \delta)$ counting
in \eqref{eq:3.6} only the cases $h^* \in \mathcal H$:
\begin{align}
\label{eq:4.5}
S_2(h, \delta) &\geq \frac{(1 + O(\delta))h^k}{N h^k(\log R)^{k + 2\ell}} \cdot \frac{N(\log R)^{k + 2\ell + 1}k}{(k + 2\ell + 1)!} {2(\ell + 1)\choose \ell + 1}\\
& \quad - \frac{(1 + O(\delta))h^k}{N h^k(\log R)^{k + 2\ell}} \cdot \frac{N(\log R)^{k + 2\ell} \log N}{(k + 2\ell)!} {2\ell \choose \ell}\nonumber\\
&= \frac{1 + O(\delta)}{(k + 2\ell)!} {2\ell\choose \ell} \left(\frac{k}{k + 2\ell + 1} \cdot
\frac{2(2\ell + 1)}{\ell + 1} \log R - \log N\right).\nonumber
\end{align}

This implies with the choice $R = N^{(\vartheta - \varepsilon)/(2 + \delta)}$
\beq
\label{eq:4.6}
S_2(h, \delta) \gg_{k, \ell, \delta}\biggl(4\left(1 - \frac{2\ell + 1}{k + 2\ell + 1}\right) \left(1 - \frac1{2\ell + 2}\right) \frac{\vartheta - \varepsilon}{2 + \delta} - 1 + O(\delta)\biggr) \mathcal L.
\eeq

It is easy to check that with the choices
\beq
\label{eq:4.7}
k = 1880, \ \ \ell = 21, \ \ \vartheta = \frac{157}{300}
\eeq
we obtain
\beq
\label{eq:4.8}
S_2(h, \delta) \gg_{k, \ell, \delta, \varepsilon} \mathcal L
\eeq
if we choose additionally $\delta$ and $\varepsilon$ sufficiently small.

On the other hand we can prove an upper bound for $S_2(h, \delta)$ as a function of $h, N, R, Q(N, h), k$ and by using Lemmas 5 and 6 and Selberg's upper bound sieve, which, combined with the lower bound \eqref{eq:4.6}, will yield a lower bound for $Q(N, h)$.
We have, namely
\begin{align}
\label{eq:4.9}
&S_2(h, \delta)\\
&\leq \frac1{N h^k(\log R)^{k + 2\ell}} \sum_{\substack{n \sim N \\
\Theta(n, h) \geq (3/2)\log N}} \Theta(n, h) {\sum_{\mathcal H}}^* \Lambda_R^2(n; \mathcal H, k + \ell)
\nonumber\\
&\leq \frac1{Nh^k(\log R)^{k + 2\ell}}\biggl(\sum_{\substack{n \sim N\\
\Theta(n, h) \geq (3/2)\log N}}\!\!\!\! 1 \biggr)^{\!\! 1/2}
\left\{\sum_{n \sim N} \Theta^2(n, h) \biggl({\sum_{\mathcal H}}^* \!\! \Lambda_R^2 (n; \mathcal H, k \!+\! \ell)\!\biggr)^{\!2}\right\}^{\!\! 1/2} \nonumber\\
&:= \frac{(Q(N, h))^{1/2}}{Nh^k(\log R)^{k + 2\ell}} I^{1/2}
\nonumber
\end{align}
where
\beq
\label{eq:4.10}
I = \!\!\! \sum_{1 \leq h', h''\leq h} \sum_{\substack{\mathcal H\subseteq [1, h], \mathcal H = k\\
\mathcal H_i \text{ admissible}\\
i = 1, 2}} \sum_{\substack{n \sim N\\
\bigl(\!P_{\mathcal H_1\cup \mathcal H_2}(n), P(R^\delta)\!\bigr) = 1}}
\!\!\!\!\!\!\!\!\!\!\!\!\! \theta(n + h') \theta(n + h'') \Lambda_R^2(n; \mathcal H_1, \ell) \Lambda_R^2(n; \mathcal H_2, \ell).
\eeq
The first important observation is that we have for any $\mathcal H$ in the above sum
\beq
\label{eq:4.11}
\Lambda_R(n; \mathcal H, k + \ell) \ll_{k, \ell, \delta} (\log R)^{k + \ell},
\eeq
since all prime divisors of $P(\mathcal H_i, n)$ are at least $R^\delta$ and $R\geq N^{c_1}$, due to the notation ${\sum\limits_{\mathcal H}}^*$ in \eqref{eq:4.4}.
This implies with the notation $\mathcal H_0 = \{h'\} \cup \{h''\} \cup \mathcal H_1 \cup \mathcal H_2$
\beq
\label{eq:4.12}
I \ll_{k, \ell, \delta} \mathcal L^2(\log R)^{4(k + \ell)} \sum_{1 \leq h', h'' \leq h} \sum_{\substack{
\mathcal H_i \subseteq [1, h]|\mathcal H_i| = k\\
\mathcal H_i \text{ admissible}\\
i = 1,2}} \sum_{\substack{n \sim N\\
\bigl(P_{\mathcal H_0}(n), P(R^\delta)\bigr) = 1}} 1.
\eeq
For a given $\mathcal H_0 \subseteq \{1, \ldots, h\}$ we have $|\mathcal H_0| = k + r$, $0 \leq r \leq k + 2$ and the total number $D(h, r)$ of quadruples $h', h'', \mathcal H_1, \mathcal H_2$ corresponding to $\mathcal H_0$ is $O_{k,r}(1) = O_k(1)$, so
\beq
\label{eq:4.13}
I \ll_{k, \ell, \delta} \mathcal L^2(\log R)^{4(k + \ell)} \sum_{r = 0}^{k + 2} \sum_{|\mathcal H_0| = k + r} \sum_{\substack{n \sim N\\
\bigl(P_{\mathcal H_0}(n) P(R^\delta)\bigr) = 1}} 1.
\eeq

For any given $\mathcal H_0$ we can apply Selberg's upper bound sieve (Theorem 5.1 of \cite{HR1974}) to obtain
\beq
\label{eq:4.14}
\sum_{\substack{n \sim N\\
(P_{\mathcal H}(n), P(R^\delta)) = 1}} 1 \leq \frac{(1 + o(1))N |\mathcal H|! \mathfrak S(\mathcal H)}{(\log R^\delta)^{|\mathcal H|}}.
\eeq

For a given $r$ summation over all $\mathcal H_0$ with $|\mathcal H_0| = k + r$ this gives by \eqref{eq:1.15}
\begin{align}
\label{eq:4.15}
\sum_{\substack{n \sim N, |\mathcal H_0| = k + r\\
\bigl(P_{\mathcal H_0}(n), P(R^\delta)\bigr) = 1}} \!\!\! 1
&\leq \sum_{|\mathcal H_0| = k + r} \frac{(1 + o(1))(k + r)! \mathfrak S(\mathcal H_0)}{(\delta \log R)^{k + r}} N\\
&\ll_{k, \delta} \sum_{|\mathcal H_0| = k + r} \frac{\mathfrak S(\mathcal H_0)}{(\log R)^{k + r}} N
\nonumber\\
&\ll \sum_{r = 0}^{k + 2} \left(\frac{h}{\log R}\right)^{k + r} \ll \left(\frac{h}{\log R}\right)^k N, \nonumber
\end{align}
due to $h \ll \log N \ll \log R$.
Hence, by \eqref{eq:4.13}--\eqref{eq:4.15} we have
\beq
\label{eq:4.16}
I^{1/2} \ll_{k, \ell, \delta} \mathcal L h^{k/2}(\log R)^{3k/2 + 2\ell} N^{1/2}.
\eeq

Together with \eqref{eq:4.6} and \eqref{eq:4.9} this implies
\beq
\label{eq:4.17}
\mathcal L \ll_{k, \ell, \delta} S_2(h, \delta) \ll_{k, \ell, \delta} \mathcal L Q(N, h)^{1/2} \left(\frac{\log R}{h}\right)^{k/2} N^{-1/2},
\eeq
consequently, by \eqref{eq:4.2},
\beq
\label{eq:4.18}
h \sum_{\substack{N < p_j \leq 2N\\
p_{j + 1} - p_j \leq h}}\!\!\! 1 \gg_{k, \ell, \delta} N\left(\frac{h}{\log R}\right)^k + O\bigl(Ne^{-c\sqrt{\mathcal L}}\bigr).
\eeq
From this we obtain by $h/\log R \geq h/\log N \geq 1/\mathcal L$
\beq
\label{eq:4.19}
\qquad \qquad \sum_{\substack{N < p_j \leq 2N\\ p_{j + 1} - p_j \leq h}} 1 \gg_{k, \ell, \delta} \frac{\eta^k N}{h} = \eta^{k - 1} \frac{N}{\mathcal L} \sim \eta^{k -1}\pi(N). \qquad \qquad \square
\eeq

\section*{Acknowledgement} This research was performed while the
first author was a guest researcher at the HUN-REN Alfr\'{e}d R\'{e}nyi
Institute of Mathematics in the frame of the MTA Distinguished Guest
Fellow Scientist Program 2024 of the Hungarian Academy of Sciences, and was supported by the Hungarian National Research Developments and Innovation Office, NKFIH, Excellence 151341.

\newpage

{\small
\noindent
\'Akos Magyar\\
HUN-REN Alfr\'ed R\'enyi Institute of Mathematics \\
Budapest, Re\'altanoda u. 13--15,
H-1053 Hungary\\
and\\
The University of Georgia\\
Athens, GA 30602, United States\\
e-mail: amagyar@uga.edu\\
\\
J\'anos Pintz\\
HUN-REN Alfr\'ed R\'enyi Institute of Mathematics \\
Budapest, Re\'altanoda u. 13--15,
H-1053 Hungary\\
e-mail: pintz@renyi.hu}

\end{document}